\documentclass[a4paper, 11pt]{amsart}   
\usepackage{mathptmx, amssymb,amscd,latexsym, eulervm}
\usepackage{changepage} 
\usepackage{amsmath}
\usepackage{amsthm}
\usepackage{mathdots}
\usepackage[pagebackref,colorlinks=true,linkcolor=blue,urlcolor=blue]{hyperref}
\usepackage{color}
\usepackage[onehalfspacing]{setspace}
\usepackage{tabularx}
\usepackage{amsfonts}
\usepackage{paralist}
\usepackage{aliascnt}
\usepackage[initials, lite]{amsrefs}
\usepackage{amscd}
\usepackage{blkarray}
\usepackage{mathbbol}
\usepackage{setspace}
\usepackage[inner=2.4cm,outer=2.4cm, bottom=3.2cm]{geometry}
\usepackage{tikz, tikz-cd}
\usepackage{calligra,mathrsfs}
\usepackage{comment}

\usepackage{tikz}
\usetikzlibrary{matrix}
\usetikzlibrary{arrows,calc}
\usetikzlibrary{decorations.markings,intersections,positioning,calc}

  \tikzset{mylabel/.style  args={at #1 #2  with #3}{
    postaction={decorate,
    decoration={
      markings,
      mark= at position #1
      with  \node [#2] {#3};
 } } } }
\allowdisplaybreaks

\BibSpec{collection.article}{%
	+{}  {\PrintAuthors}                {author}
	+{,} { \textit}                     {title}
	+{.} { }                            {part}
	+{:} { \textit}                     {subtitle}
	+{,} { \PrintContributions}         {contribution}
	+{,} { \PrintConference}            {conference}
	+{}  {\PrintBook}                   {book}
	+{,} { }                            {booktitle}
	+{,} { }                            {series}
	+{, vol.} { }                            {volume}
	+{,} { }                            {publisher}
	+{,} { \PrintDateB}                 {date}
	+{,} { pp.~}                        {pages}
	+{,} { }                            {status}
	+{,} { \PrintDOI}                   {doi}
	+{,} { available at \eprint}        {eprint}
	+{}  { \parenthesize}               {language}
	+{}  { \PrintTranslation}           {translation}
	+{;} { \PrintReprint}               {reprint}
	+{.} { }                            {note}
	+{.} {}                             {transition}
	+{}  {\SentenceSpace \PrintReviews} {review}
}
\AtBeginDocument{%
	\def\MR#1{}
}

\makeatletter
\@namedef{subjclassname@2020}{%
	\textup{2020} Mathematics Subject Classification}
\makeatother

\newcommand{\tr}{\operatorname{tr}}

\def\opn#1#2{\def#1{\operatorname{#2}}}
\opn\Cl{Cl} \opn\deg{deg} \opn\Stab{Stab} \opn\aff{aff} \opn\div{div}
\opn\cone{cone} \opn\End{End} \opn\mod{mod}  \opn\pdim{pdim} \opn\diag{diag} \opn\vert{vert} \opn\m{m} \opn\V{V}
\opn\Cone{Cone} \opn\Pyr{Pyr} \opn\max{max} \opn\min{min} \opn\int{int} \opn\rev{rev} \opn\ker{ker} \opn\lat{lat} \opn\pull{pull}
\opn\cok{coker} \opn\ant{ant}
\opn\inte{int}

\newcommand{\kk}{\mathbb{k}}

\newcommand{\NN}{\normalfont\mathbb{N}}
\newcommand{\ZZ}{\normalfont\mathbb{Z}}

\newcommand{\mm}{{\normalfont\mathfrak{m}}}
\newcommand{\MM}{{\normalfont\mathfrak{m}}}

\newcommand{\pp}{{\normalfont\mathfrak{p}}}

\newcommand{\nn}{{\normalfont\mathfrak{n}}}

\newcommand{\depth}{\normalfont\text{depth}}

\newcommand{\Ann}{\normalfont\text{Ann}}

\newcommand{\Ass}{{\normalfont\text{Ass}}}

\newcommand{\Hom}{\normalfont\text{Hom}}

\newcommand{\indeg}{\normalfont\text{indeg}}

\newcommand{\Hilb}{{\normalfont\text{Hilb}}}
\newcommand{\Spec}{\normalfont\text{Spec}}

\def\f0{\mathbf{0}}

\def\1{\mathbf{1}}

\newtheorem{theorem}{Theorem}[section]

\newaliascnt{headcor}{headthm}

\aliascntresetthe{headcor}

\newaliascnt{headconj}{headthm}

\aliascntresetthe{headconj}

\newaliascnt{corollary}{theorem}
\newtheorem{corollary}[corollary]{Corollary}
\aliascntresetthe{corollary}

\newaliascnt{claim}{theorem}

\aliascntresetthe{claim}

\newaliascnt{lemma}{theorem}
\newtheorem{lemma}[lemma]{Lemma}
\aliascntresetthe{lemma}

\newaliascnt{conjecture}{theorem}

\aliascntresetthe{conjecture}

\newaliascnt{proposition}{theorem}

\aliascntresetthe{proposition}

\theoremstyle{definition}
\newaliascnt{definition}{theorem}
\newtheorem{definition}[definition]{Definition}
\aliascntresetthe{definition}

\newaliascnt{notation}{theorem}
\newtheorem{notation}[notation]{Notation}
\aliascntresetthe{notation}

\newaliascnt{example}{theorem}

\aliascntresetthe{example}

\newaliascnt{examples}{theorem}

\aliascntresetthe{examples}

\newaliascnt{remark}{theorem}
\newtheorem{remark}[remark]{Remark}
\aliascntresetthe{remark}

\newaliascnt{question}{theorem}

\aliascntresetthe{question}

\newaliascnt{questions}{theorem}

\aliascntresetthe{questions}

\newaliascnt{problem}{theorem}

\aliascntresetthe{problem}

\newaliascnt{construction}{theorem}

\aliascntresetthe{construction}

\newaliascnt{setup}{theorem}
\newtheorem{setup}[setup]{Setup}
\aliascntresetthe{setup}

\newaliascnt{algorithm}{theorem}

\aliascntresetthe{algorithm}

\newaliascnt{observation}{theorem}

\aliascntresetthe{observation}

\newaliascnt{defprop}{theorem}

\aliascntresetthe{defprop}

\def\equationautorefname~#1\null{(#1)\null}
\def\sectionautorefname~#1\null{Section #1\null}
\def\subsectionautorefname~#1\null{\S #1\null}

\title{When do pseudo-Gorenstein rings become Gorenstein?}

\author[S. Miyashita]{Sora Miyashita}
\address[S. Miyashita]{Department of Pure And Applied Mathematics, Graduate School Of Information Science And Technology, Osaka University, Suita, Osaka 565-0871, Japan}
\email{u804642k@ecs.osaka-u.ac.jp}

\date{\today}

\begin{document}

	\maketitle

\begin{abstract}
We discuss the relationship between the trace ideal of the canonical module and pseudo-Gorensteinness.
In particular, under certain mild assumptions, we show that every positively graded domain that is both pseudo-Gorenstein and nearly Gorenstein is Gorenstein.
As an application, we clarify the relationships among nearly Gorensteinness, almost Gorensteinness, and levelness---notions that generalize Gorensteinness---in the context of standard graded domains.
Moreover, we give a method for constructing quasi-Gorenstein rings by taking a Veronese subalgebra of certain Noetherian graded rings.
\end{abstract}

\section{Introduction}
Throughout this introduction,
let $R = \bigoplus_{i \ge 0} R_i$ be a Cohen--Macaulay positively graded ring with $R_0$ a field.
Let $\MM_R=\bigoplus_{i>0} R_i$ be the graded maximal ideal of $R$.
Furthermore, let $\omega_R = \bigoplus_{i \in \ZZ} [\omega_R]_i$ be the graded canonical module of $R$
and let $a_R=-\min\{i \in \ZZ : [\omega_R]_i \neq 0\}$ be the {\it $a$-invariant} of $R$.
\begin{definition}
Based on \cite{ene2015pseudo},
we say that $R$ is {\it pseudo-Gorenstein}~if 
$\dim_{R_0}([\omega_R]_{-a_R}) = 1$.
\end{definition}
The pseudo-Gorenstein property in certain classes of Cohen--Macaulay standard graded rings has been actively studied (see, e.g., \cite{hatasa2024pseudo,ene2015pseudo,rinaldo2024level,rinaldo2023level}).
Every Gorenstein graded ring is pseudo-Gorenstein, but not vice versa.
\emph{Under what non-trivial conditions is a pseudo-Gorenstein ring Gorenstein?}

{\it Level rings}, introduced by Stanley,
are among the earliest generalizations of Gorenstein rings.
A ring $R$ is called {\it level}~(\cite{stanley2007combinatorics}) if the minimal generators of $\omega_R$ all have the same degree.
Note that if $R$ is level and pseudo-Gorenstein, then it is Gorenstein.
The third generalized Gorenstein property of interest in this paper is the {\it nearly Gorenstein property}, introduced in \cite{herzog2019trace}.
Let $\tr_R(\omega_R)$ denote the ideal generated by the images of all graded $R$-linear maps from $\omega_R$ to $R$.  
$R$ is called {\it nearly Gorenstein}~(\cite{herzog2019trace}) if $\tr_R(\omega_R) \supseteq \MM_R$.  

To state the main results, we first need to introduce some terminology and notation. 
We say that \( R \) is \emph{generically Gorenstein} if \( R_{\mathfrak{p}} \) is Gorenstein for every minimal prime ideal \( \mathfrak{p} \) of \( R \).  
An element \( x \in \omega_R \) is said to be \emph{torsion-free} if \( rx = 0 \) implies \( r = 0 \) for all \( r \in R \).  
We define \( \indeg_R(\mathfrak{m}_R) := \min \{ i \in \mathbb{Z} : [\mathfrak{m}_R]_i \neq 0 \} \).  
Our main result is stated as follows.

\begin{theorem}
[{\autoref{mainTHM}}]
\label{XXX}
Assume that $R$ is generically Gorenstein, that $\omega_R$ contains a torsion-free homogeneous element of degree $-a_R$, and that $\tr_R(\omega_R)$ contains an $R$-regular sequence $\theta_1,\theta_2 \in R_{\indeg_R(\MM_R)}$.
If $R$ is pseudo-Gorenstein, then it is Gorenstein.
\end{theorem}

For instance, if $R$ is a nearly Gorenstein standard graded domain of dimension at least two, and $R_0$ is an infinite field, then $R$ satisfies all the assumptions of \autoref{XXX}.
This shows that, even among nearly Gorenstein rings, pseudo-Gorensteinness often implies Gorensteinness.

Since the introduction of nearly Gorenstein rings, the study of these rings and the trace ideals has garnered significant attention, becoming an active area of research in recent years~(\cite{hall2023nearly,dao2020trace,miyashita2024linear,herzog2019trace,ficarra2024canonical,ficarra2024canonical!,caminata2021nearly,kumashiro2023nearly,bagherpoor2023trace,lu2024chain,lyle2024annihilators,miyazaki2024non,jafari2024nearly,miyashita2024canonical,celikbas2025full,kimura2025trace}).
In addition to the aforementioned levelness and nearly Gorensteinness, one major generalization of Gorensteinness is \emph{almost Gorensteinness}
(\cite{goto2013almost,goto2015almost,barucci1997one}).
The relationships among these three properties have been actively investigated across various classes of rings~(\cite{herzog2019trace,miyashita2024comparing,moscariello2021nearly,miyashita2024canonical,hall2023nearly,higashitani2022levelness}).
In the one-dimensional case, it has been established that every almost Gorenstein ring is necessarily nearly Gorenstein~(\cite[Proposition 6.1]{herzog2019trace}). Furthermore, if $R$ has minimal multiplicity and positive Krull dimension, and if $R_0$ is an infinite field, then $R$ is almost Gorenstein if it is nearly Gorenstein~\cite[Theorem 6.6]{herzog2019trace}.
However, in the case of higher-dimensional rings or those without minimal multiplicity, the relationships among these properties are still not well understood. 

As mentioned before, level rings are closely related to pseudo-Gorenstein rings. Similarly, it is known that almost Gorenstein rings are also strongly related to pseudo-Gorenstein rings in the case of standard graded domains~(\cite[Theorem 4.7]{higashitani2016almost}).
A new perspective on the relationship between nearly Gorenstein rings and pseudo-Gorenstein rings is provided by \autoref{XXX}.
By focusing on the pseudo-Gorenstein property, we clarify how levelness, nearly Gorensteinness, and almost Gorensteinness relate to one another in domains of higher dimension.
Let \( r(R) \) denote the \emph{Cohen--Macaulay type} of \( R \).

\begin{theorem}[{\autoref{YYYY} and \autoref{Thm:MainC}}]\label{YYY}
Suppose that \(R\) is a nearly Gorenstein graded domain with \(\dim(R) \ge 2\), such that the ideal generated by \(R_{\indeg_R(\mm_R)}\) is \(\mm_R\)-primary.
Then the following hold:
\begin{itemize}
\item[\rm (1)] $R$ is pseudo-Gorenstein if and only if it is Gorenstein;
\item[\rm (2)] If $r(R)=2$, then $R$ is level;
\item[\rm (3)] Moreover, assume that \( R \) is a standard graded almost Gorenstein ring, and \( R_0\) is an algebraically closed field of characteristic zero.  
Then \( R \) must be either Gorenstein or have minimal multiplicity.
In particular, \( R \) is level in either case.
\end{itemize}
\end{theorem}


Indeed, \autoref{XXX} extends to a more general setting of Noetherian rings that are not necessarily Cohen--Macaulay (see \autoref{Thm:NonCM2}).  
Finally, we present a result obtained by applying \autoref{Thm:NonCM2} to Veronese subalgebras.
If $R$ is a standard graded Gorenstein ring of positive dimension, then any Veronese subalgebra of $R$ is nearly Gorenstein (see \cite[Corollary 4.7]{herzog2019trace}).
Note that one cannot construct a Gorenstein ring from the Veronese subalgebra of a non-Gorenstein Cohen–Macaulay standard graded ring of dimension at least two
(\cite[Theorem (3.2.1)]{goto1978graded}).
However, we show that by considering the analogous properties of pseudo-Gorenstein (resp. nearly Gorenstein) rings for a Noetherian ring, one can construct quasi-Gorenstein rings~(see \autoref{def:qG}) from certain Veronese subalgebras.
In the following, we do not assume that \( A \) is Cohen--Macaulay. Since the canonical module \(\omega_A\) and the \(a\)-invariant \(a_A\) of \( A \) can be defined similarly,
we use the same notation as in the Cohen--Macaulay case.

\begin{theorem}[{\autoref{OKOKOK}}]
\label{ZZZ}
Let \( A=\bigoplus_{i \ge 0} A_i\) be a Noetherian standard graded domain with $A_0$ a field and let $\mm_A=\bigoplus_{i>0} A_i$.
Let \( k \) be a positive integer
and let $A^{(k)}=\bigoplus_{i \ge 0} A_{ik}$.
Assume that \( a_A \in k\mathbb{Z} \) and that \( \depth(A^{(k)}) \geq 2 \).  
Suppose further that
\[
  \tr_A(\omega_A) \supseteq \MM_A
  \quad \text{and} \quad
  \dim_{A_0}([\omega_A]_{-a_A}) = 1.
\]
Then \( A^{(k)} \) is quasi-Gorenstein.
In particular, if \( A^{(k)} \) is Cohen--Macaulay, then it is Gorenstein.
\end{theorem}

\subsection*{Outline}
In \autoref{sectA}, we prepare some foundational facts and definitions for later discussions, collecting basic definitions and properties about generalizations of the Gorenstein property.
In \autoref{sectB}, we prove the general statement \autoref{Thm:General} for $R$-modules that are not necessarily finitely generated over a graded ring $R$. As a consequence, we establish \autoref{XXX}.
In \autoref{sectC}, as an application of \autoref{XXX}, we discuss various corollaries derived from it. First, we clarify the relationships among nearly Gorensteinness, almost Gorensteinness and levelness in the setting of standard graded domains~(\autoref{YYY}). Furthermore, we explore rings that are non-Cohen--Macaulay but satisfy an analogue of the nearly Gorenstein property, demonstrating that their Veronese subalgebras yield quasi-Gorenstein rings~(\autoref{ZZZ}).

\begin{setup}\label{setup1}
Throughout this paper,
we denote the set of non-negative integers by $\NN$.
Let $R=\bigoplus_{i \in \NN} R_i$ be a positively graded ring. Unless otherwise stated, we assume that $R_0$ is a field. Hence, $R$ has the unique graded maximal ideal given by $\MM_R:=\bigoplus_{i \in \NN \setminus \{0\}} R_i$.
Let $M$ and $N$ be graded $R$-modules
and let $k \in \ZZ$.
Let $M(k)$ denote the graded $R$-module having the same underlying $R$-module structure as $M$,
where $[M(k)]_n = [M]_{n+k}$ for all $n \in \ZZ$.
Moreover, let ${}^*\Hom_R(M, N)$ denote the graded $R$-module consisting of graded homomorphisms from $M$ to $N$.
Let \( \mathcal{S} \subseteq R \) be the multiplicative set of homogeneous non-zero divisors. 
The \emph{graded total quotient ring} \( {}^*Q(R) \) is defined as the localization
${}^*Q(R) := \mathcal{S}^{-1}R$.
Suppose that $R$ is a Noetherian ring. 
Let $\omega_R=\bigoplus_{i \in \ZZ} [\omega_R]_i$ denote the graded canonical module of $R$ (see \autoref{defcanon}). Let $a_R$ denote the {\it $a$-invariant} of $R$, that is, 
\[
a_R := - \min \{ i \in \mathbb{Z} : [\omega_R]_i \neq 0 \}.
\]
\end{setup}

\section{Preliminaries}\label{sectA}
The purpose of this section is to lay the groundwork for the proof of our main result. Throughout, we work under \autoref{setup1} in all subsections. Furthermore, except in the final subsection \autoref{lastsubsec}, we assume that $R$ is a Noetherian ring.




\subsection{Canonical modules over Noetherian graded rings}
Let us recall the definition of the canonical module over a positively graded Noetherian ring.
Throughout this subsection, we retain \autoref{setup1}, and we further assume that \( R \) is Noetherian.

\begin{definition}{\cite[Definition (2.1.2)]{goto1978graded}}\label{defcanon}
Let $d=\dim (R)$. The finitely generated graded $R$-module
\[
\omega_R:={}^*\Hom_R({}^*\mathrm{H}_{\MM_R}^{d}(R), E_R)
\]
is called the {\it graded canonical module}, where ${}^*\mathrm{H}_{\MM_R}^{d}(R)$ denotes the graded $d$-th local cohomology and $E_R$ denotes the graded injective envelope of $R_0$.
\end{definition}

\begin{definition}\label{def:qG}
A ring \( R \) is said to be \emph{quasi-Gorenstein} if \( \omega_R \) is isomorphic to \( R \) as a graded \( R \)-module (up to degree shift).
\end{definition}

\begin{remark}
A ring $R$ is Gorenstein if and only if it is Cohen–Macaulay and quasi-Gorenstein.
\end{remark}

\begin{definition}\label{def:pseudo}
Assume that $R$ is Cohen--Macaulay.
Following the terminology of \cite{ene2015pseudo},
we say that $R$ is {\it pseudo-Gorenstein} if
$\dim_{R_0}([\omega_R]_{-a_R})=1$.
Moreover,
$R$ is called {\it level}~(\cite{stanley2007combinatorics}) if
$\omega_R$ is generated in a single degree as a graded $R$-module.
\end{definition}

\begin{remark}
Being pseudo-Gorenstein does not imply being Gorenstein in general.
For example, all numerical semigroup rings (see, e.g., \cite[p.~178]{bruns1998cohen} for the definition) are pseudo-Gorenstein in the sense of \autoref{def:pseudo}, but infinitely many of them are not Gorenstein~(see \cite[Theorem 4.4.8]{bruns1998cohen}).  
On the other hand, every graded Cohen–Macaulay ring that is both level and pseudo-Gorenstein is Gorenstein by definition.
\end{remark}

\subsection{Semi-standard graded rings, Hilbert series, and the graded almost Gorenstein property}
We now state an important result (\autoref{WOW}) concerning the almost Gorenstein property and pseudo-Gorensteinness defined in the previous section.
Throughout this subsection, we retain \autoref{setup1}, and we further assume that \( R \) is Noetherian.

\begin{definition}
If $R =R_0[R_1]$, that is, if $R$ is generated by $R_1$ as an $R_0$-algebra, then $R$ is called {\it standard graded}.
If $R$ is finitely generated as an $R_0[R_1]$-module, then $R$ is called {\it semi-standard graded}.
\end{definition}

\begin{remark}
The concept of semi-standard graded rings, which generalizes standard graded rings, arises naturally in this context from the perspective of combinatorial commutative algebra.
For instance, the \emph{Ehrhart rings} of lattice polytopes and the \emph{face rings} of simplicial posets (see \cite{stanley1991f}) are typical examples of semi-standard graded rings.
\end{remark}

Recall that if $R$ is semi-standard graded, then its Hilbert series $\Hilb(R,t) := \sum_{i=0}^\infty \dim_{R_0}(R_i)t^i$ is of the following form:
\[
\Hilb(R,t) = \frac{h_0 + h_1 t + \cdots + h_s t^s}{(1-t)^{\dim R}},
\]
where $h_i \in \mathbb{Z}$ for all $0 \le i \le s$, $h_s \ne 0$, and $\sum_{i=0}^s h_i \ne 0$.
The integer \( s \), known as the {\it socle degree} of \( R \), is denoted by \( s(R) \).  
The sequence \(h(R)=(h_0, h_1, \ldots, h_{s(R)}) \) is called the {\it \( h \)-vector} of \( R \).

\begin{remark}
Assume that $R$ is semi-standard graded, and let $h(R)=(h_0, h_1, \ldots, h_{s(R)})$ be the $h$-vector of $R$.
We always have $h_0=1$ because $h_0=\Hilb(R,0)=\dim_{R_0}(R_0)=1$.
Moreover, if $R$ is Cohen--Macaulay,  its $h$-vector satisfies $h_{s(R)}=\dim_{R_0}([\omega_R]_{-a_R})$
and $h_i \geqq 0$ for every $0 \le i \le s(R)$.
Therefore, in this case, \( R \) is pseudo-Gorenstein if and only if \( h_{s(R)} = 1 \).
For further information on the $h$-vectors of Noetherian graded rings, see \cite{stanley1978hilbert} and \cite{stanley2007combinatorics}.
\end{remark}


\begin{definition}
Assume that $R$ is Cohen--Macaulay.
$R$ is called {\it almost Gorenstein}~(\cite[Section 10]{goto2015almost}) if there exists a graded $R$-monomorphism $\phi: R \hookrightarrow \omega_R(-{a_R})$ of degree 0 such that $C:=\cok(\phi)$ is either the zero module or an Ulrich $R$-module,
i.e. $C$ is Cohen--Macaulay and $\mu(C) = e(C)$.
Here, $\mu(C)$ (resp. $e(C)$) denotes the minimal number of generators of $C$ (resp. the multiplicity of $C$ with respect to $\MM_R$).
\end{definition}


\begin{theorem}[{\cite[Theorem 4.7]{higashitani2016almost}}]\label{WOW}
Suppose that $R$ is a Cohen--Macaulay standard graded domain
and $R_0$ is an algebraically closed field with characteristic 0.
If $R$ is almost Gorenstein,
then it is pseudo-Gorenstein,
or it has a minimal multiplicity.
\begin{proof}
Let $h(R)=(h_0,\cdots,h_{s(R)})$ be the $h$-vector of $R$.
When $R_0$ is an infinite field, observe that $R$ is pseudo-Gorenstein if and only if $h_{s(R)} = 1$, and that $R$ has minimal multiplicity if and only if $s(R) \le 1$.  
In light of this, the statement is merely a reformulation of \cite[Theorem 4.7]{higashitani2016almost}.
\end{proof}
\end{theorem}

\subsection{Trace ideals over graded rings}\label{lastsubsec}
Throughout this subsection, we retain \autoref{setup1}, except that \( R \) is not necessarily assumed to be Noetherian.

\begin{definition}\label{def:GRADEDcanon}
Let $M$ be a (not necessarily finitely generated) graded $R$-module.
The sum of all images of graded homomorphisms $\phi \in {}^*\Hom_R(M,R)$ is called the {\it trace ideal} of $M$:
\[
\tr_R(M):=\sum_{\phi \in {}^*\Hom_R(M,R)}\phi(M).
\] 
\end{definition}

\begin{remark}[{\cite[Remark 2.2]{kumashiromiyashita2025}}]
Let $M$ be a  graded $R$-module.
If $R$ is Noetherian,
then we have
$$\tr_R(M)=\sum_{\phi \in \Hom_R(M,R)} \phi(M).$$
\end{remark}

\begin{remark}
Let $I$ be a graded ideal of $R$ containing
a non-zero divisor of $R$.
Set
$$\displaystyle I^{-1}:=\{x \in {}^*Q(R):xI\subseteq R\}.$$
Then we have $\tr_R(I)=I \cdot I^{-1}$
(see the proof of \cite[Lemma 1.1]{herzog2019trace}).
\end{remark}

\begin{remark}\label{rem:gradedaoyamagoto}
Assume that $R$ is Noetherian.
Then $R$ is quasi-Gorenstein if and only if $\tr_R(\omega_R)=R$
~(see \cite[Remark 2.9~(4)]{kumashiromiyashita2025}).
Moreover,
\(\tr_R(\omega_R)\) describes the non-quasi-Gorenstein locus if \(R\) is {\it unmixed},
that is,
$\dim(R)=\dim(R/\pp)$ for all $\pp \in \Ass(R)$
~(see \cite[Remark 2.9~(2), (5)]{kumashiromiyashita2025}).
It is known that $R$ is unmixed if and only if $\Ann_R(\omega_R)=(0)$~(see \cite[Remark 2.9~(2), (3)]{kumashiromiyashita2025}).
\end{remark}

\begin{definition}[{\cite[Definition 2.2]{herzog2019trace}}]
Assume that $R$ is Cohen--Macaulay.
$R$ is called \textit{nearly Gorenstein} if
$\tr_R(\omega_R) \supseteq {\MM_R}$.
\end{definition}

\begin{definition}[{\cite[Definition 1.1]{miyashita2024linear}}]
Assume that \(R\) is a semi-standard graded Cohen--Macaulay ring.  
We say that {\it \(R\) satisfies \((\natural)\)} if \(\sqrt{[\tr_R(\omega_R)]_1 R} \supseteq \MM_R\),  
i.e., the radical of the ideal generated by the degree-one part of \(\tr_R(\omega_R)\) contains \(\mm_R\).
\end{definition}

\begin{remark}\label{aaaaaaaaaaaaaaaaaaaaaaaaaaaaaaaaaaaaaaaaaaaaaaa}
Assume that $R$ is a semi-standard graded Cohen--Macaulay ring,
and that $R_0$ is an infinite field.
If
$R$ is nearly Gorenstein,
then it satisfies $(\natural)$.
\begin{proof}
It follows from \cite[Proposition 1.5.12]{bruns1998cohen}.
\end{proof}
\end{remark}

The result below is known when \( R \) is Cohen--Macaulay or a local ring.  
However, in the case where \( R \) is positively graded and not Cohen--Macaulay, the author is not aware of an appropriate reference.  
Therefore,
we provide a proof for completeness.

\begin{remark}\label{HelpMeFromMugentai}
Assume that $R$ is Noetherian.
Let $\kk$ be a field and set $A = R \otimes_{R_0} \kk$.
Then the following hold:
\begin{itemize}
\item[\rm (1)] $\omega_A \cong \omega_R \otimes_{R_0} \kk$ as a graded $A$-module;
\item[\rm (2)] $\tr_A(\omega_A)=\tr_R(\omega_R)A$;
\end{itemize}
\begin{proof}
(1): Let \( \omega_{R_{\mathfrak{m}_R}} \) and \( \omega_{A_{\mathfrak{m}_A}} \) denote the canonical modules of the Noetherian local rings \( R_{\mathfrak{m}_R} \) and \( A_{\mathfrak{m}_A} \), respectively (see \cite[Definition~1.1]{aoyama1983some}).
Then, we obtain \( \omega_{A_{\mathfrak{m}_A}} \cong \omega_{R_{\mathfrak{m}_R}} \otimes_{R_0} \mathbb{k} \) by \cite[Theorem~4.1]{aoyama1985endomorphism}.
On the other hand,
by \cite[Remark~2.9~(1)]{kumashiromiyashita2025}, we have
$$(\omega_A)_{\mathfrak{m}_A} \cong \omega_{A_{\mathfrak{m}_A}},
\;\;\;
\omega_{A_{\mathfrak{m}_A}} \cong (\omega_A)_{\mathfrak{m}_A},
\;\;\;\text{and}\;\;\;
\omega_{R_{\mathfrak{m}_R}} \otimes_{R_0} \mathbb{k} \cong (\omega_R)_{\mathfrak{m}_R} \otimes_{R_0} \mathbb{k} \cong (\omega_R \otimes_{R_0} \mathbb{k})_{\mathfrak{m}_A}$$ as \( A_{\mathfrak{m}_A} \)-modules.
Therefore, we conclude that \( (\omega_A)_{\mathfrak{m}_A} \cong (\omega_R \otimes_{R_0} \mathbb{k})_{\mathfrak{m}_A} \) as an \( A_{\mathfrak{m}_A} \)-module, and hence, by \cite[Corollary~3.9]{hashimoto2023indecomposability}, we obtain \( \omega_A \cong \omega_R \otimes_{R_0} \mathbb{k} \) as a graded \( A \)-module.

(2):
It follows from (1) and \cite[Lemma 1.5~(iii)]{herzog2019trace}.
\end{proof}
\end{remark}

At the end of this section, we recall the definition of the Veronese subalgebra and present a lemma concerning its relation to the trace ideal of the canonical module.

\begin{definition}
Let \( k > 0 \) be a positive integer. For \( R = \bigoplus_{i \in \NN} R_i \), the \textit{\( k \)-th Veronese subalgebra} of $R$ is defined by
$R^{(k)} := \bigoplus_{i \in \NN} R_{ik}$.
This may seem tautological, we emphasize that the grading on \( R^{(k)} \) is given by \( [R^{(k)}]_i = R_{ik} \) for each \( i \in \NN \).  
Similarly, for a graded \( R \)-module \( M = \bigoplus_{i \in \ZZ} M_i \), we define its \textit{\( k \)-th Veronese submodule} of $M$ by
$M^{(k)} := \bigoplus_{i \in \ZZ} M_{ik}$.
\end{definition}

The following is the non-Cohen–Macaulay version of \cite[Proposition 6.2]{miyashita2024linear}.

\begin{lemma}\label{lem:Yeah}
Assume that \( R \) is Noetherian and generically Gorenstein, and that there exists a non-zero divisor \( z \in R_1 \).  
Let \( I \) be a graded ideal generated by a subset of \( R_1 \), and let \( \mathfrak{n}_R \) be the ideal of \( R \) generated by \( R_1 \).  
If \( \operatorname{tr}_R(\omega_R) \supseteq I \),  
then for any integer \( k > 0 \), we have  
$\operatorname{tr}_{R^{(k)}}(\omega_{R^{(k)}}) \supseteq (\mathfrak{n}_R^{k-1} I)^{(k)}$.
\begin{proof}
By \cite[Remark 2.13]{kumashiromiyashita2025},
there exists a graded ideal \( I \) of \( R \) such that
$\omega_R \cong I$ (up to a shift).
Thus, there exists a graded fractional ideal \( J := z^{-a_R-\indeg_R(I)} I \subseteq  {}^*Q(R) \) such that \( \omega_R \cong J \) and \( -a_R = \indeg_R(J) \).
Since $\omega_{R ^{(k)}} \cong (\omega_{R}) ^{(k)}$ (see \cite[Corollary (3.1.3)]{goto1978graded}) for any $k > 0$, we have $\tr_{R^{(k)}}(\omega_{R ^{(k)}})=\tr_{R^{(k)}}(J^{(k)}) \supseteq (\nn_R^{k-1}I)^{(k)}$
by \cite[Theorem 6.1]{miyashita2024linear}.
\end{proof}
\end{lemma}

\section{The relationship between pseudo-Gorensteinness and the trace ideal of the canonical module}\label{sectB}
In this section, we prove our main result \autoref{Thm:General}.
Throughout this section, we retain \autoref{setup1}.
The following fundamental lemma plays an important role in proving our main result.

\begin{lemma}\label{lem:useful}
Let $I \subseteq R$ be a non-zero graded ideal of $R$.
Suppose that there exist homogeneous non-zero divisors $f_1, f_2 \in I$ with $f_2 \notin Rf_1$.
Let
$\theta_1,\theta_2 \in R_{\indeg_R(\MM_R)}$ be an $R$-regular sequence.
Then we have $\theta_{i} \frac{f_2}{f_1} \in Q(R) \setminus R$ for some $i=1,2$.
\end{lemma}
\begin{proof}
Suppose that $\theta_i \frac{f_2}{f_1} \in R$
for each $i=1,2$.
Then there are $r_1,r_2 \in R$ such that $\theta_1 f_2 = f_1 r_1$
and
$\theta_2 f_2 = f_1 r_2$.
Thus,
since $\theta_2 (f_1 r_1)=(f_2 \theta_1) \theta_2=(f_2 \theta_2) \theta_1=(f_1r_2) \theta_1$,
we have $f_1 r_1 \theta_2=f_1r_2 \theta_1$.
Since $f_1 \neq 0$ is a non-zero divisor of $R$,
we have $r_1 \theta_2=r_2 \theta_1$.
Moreover,
we obtain $r_1=r\theta_1$ for some $r \in R$ because $\theta_2+(\theta_1) \in R/(\theta_1)$ is a non-zero divisor of $R/(\theta_1)$.
Substituting this into $f_2 \theta_1 = f_1 r_1$,
we get
$f_2 \theta_1 = r f_1 \theta_1$.
Since $\theta_1$ is a non-zero divisor of $R$,
we have $f_2 = r f_1$.
This contradicts $f_2 \notin Rf_1$.
\end{proof}

To state the main result in a more general form, we introduce some terminology.

\begin{definition}
Let $M$ be a graded $R$-module. For $i \in \ZZ$, we say that {\it $M$ has an $R(i)$-free summand} if there exists an $R$-epimorphism $\phi \in [{}^*\Hom_R(M,R)]_i=\{ \phi \in \Hom_R(M, N) : \phi(M_k) \subseteq N_{k+i} \text{ for all } k \in \mathbb{Z} \}$.
\end{definition}

\begin{definition}
We introduce the symbols $\infty$ and $-\infty$, and define a total order with $-\infty < i < \infty$ for all $i \in \mathbb{Z}$.  
With this ordering, $\mathbb{Z} \cup \{ \infty, -\infty \}$ becomes a totally ordered set.
For a graded $R$-module $M$,
we set
\[
\indeg_R(M) :=
\begin{cases}
\infty & \text{if \;} M=(0), \\
-\infty & \text{if there exists $k \in \ZZ$ such that\;} M_i \neq (0) \text{\;for all integers $i<k$}, \\
\min\{i \in \ZZ : M_i \neq 0\} & \text{otherwise.} 
\end{cases}
\]
\end{definition}

We establish
the following, which presents \autoref{XXX} in a more general form.

\begin{theorem}\label{Thm:General}
Let $I \subseteq R$ be a graded ideal of $R$.
Suppose that there exists a non-zero divisor $f \in I_{\indeg_R(I)}$ such that $Rf \neq I$.
Let $M$ be a graded $R$-module such that $\indeg_R(M)<\infty$.
Suppose that
there exists a graded $R$-module $N$ such that
$M \cong I \oplus N$ as a graded $R$-module (up to degree shift),
$\indeg_R(I)<\indeg_R(N)$
and $[\tr_R(N)]_{\indeg_R(\MM_R)}=0$.
If $\tr_R(M)$ contains an $R$-regular sequence $\theta_1,\theta_2 \in R_{\indeg_R(\MM_R)}$ and $\dim_{R_0}(M_{\indeg_R(M)})=1$,
then $I$ has an $R(-\indeg_R(I))$-free summand.
\end{theorem}
\begin{proof}
Note that
$\indeg_R(I)=1$
because of
the graded isomorphism
$M \cong I \oplus N$,
$\indeg_R(M)=1$
and
$\indeg_R(I)<\indeg_R(N)$.
Thus $f$ is an ${R_0}$-basis of $I_{\indeg_R(I)}$.
Since $Rf \neq I$, we can choose a homogeneous element $f' \in I \setminus Rf$.
Then we have $\theta_{k} \frac{f'}{f} \notin R$ for some $k \in \{1,2\}$ by \autoref{lem:useful}.

Since ${}^{*}\operatorname{Hom}_R(I \oplus N, R) \cong {}^{*}\operatorname{Hom}_R(I, R) \oplus {}^{*}\operatorname{Hom}_R(N, R)$, it follows that $\tr_R(M) = \tr_R(I) + \tr_R(N)$.
Furthermore, since $\tr_R(I) = I \cdot I^{-1}$ by \cite[Lemma 1.1]{herzog2019trace}, it follows that $\tr_R(M) = I \cdot I^{-1} + \tr_R(N)$.
Therefore, since $\theta_{k} \in \tr_R(M)$ and $f$ is a ${R_0}$-basis of $I_{\indeg_R(I)}$, there exist
homogeneous elements
$x_1,\cdots,x_r \in I$,
$g_0, g_1,\cdots,g_r \in I^{-1}$,
$y_1,\cdots,y_s \in N$
and $\phi_1,\cdots,\phi_s \in \operatorname{Hom}_R^{*}(N, R)$ 
such that
$\deg(f)<\deg(x_i)$ for any $1 \le i \le r$, and
\begin{equation}\label{eq:Einstein2}
\theta_{k}=fg_0+\sum_{i=1}^r x_ig_i+\sum_{j=1}^s\phi_j(y_j).
\end{equation}
Moreover,
we may assume that
$\deg(fg_0)=\deg(x_ig_i)=\deg(\phi_j(y_j))=\deg(\theta_{k})=\indeg_R(\MM_R)$
for every $1\le i \le r$ and $1\le j \le s$
by \autoref{eq:Einstein2}.
Then, we have 
$\phi_j(y_j) = 0$ for any $1 \le  j \le s$.

Indeed, if $\phi_{j_0}(y_{j_0}) \neq 0$ for some $1 \le  j_0 \le s$,
we have $[\tr_R(N)]_{\indeg_R(\MM_R)} \neq 0$
because
$\phi_{j_0}(y_{j_0}) \in \tr_R(N)$ is a non-zero homogeneous element of
$\deg(\phi_{j_0}(y_{j_0}))=\deg(\theta_{k})=\indeg_R(\MM_R)$.
This contradicts to $[\tr_R(N)]_{\indeg_R(\MM_R)} = 0$.

Thus, we obtain $\displaystyle \theta_{k}=fg_0+\sum_{i=1}^r x_ig_i.$
Then we can check
$x_{i_0}g_{i_0} \neq 0$ for some $1 \le  i_0 \le r$.

Indeed, assume that $x_ig_i=0$ for any $1 \le  i \le r$.
Since
$\theta_{k} = fg_0$, $f' \in I$ and $g_0 \in I^{-1}$,
we have $\theta_{k} \frac{f'}{f}=
\frac{(fg_0) f'}{f}=
g_0 f' \in R.$
This contradicts to $\theta_{k} \frac{f'}{f} \notin R$.

Thus there exists $x_{i_0}g_{i_0} \neq 0$ for some $1 \le  i_0 \le r$.
Consequently, we can check $fg_{i_0} \neq 0$.

Indeed, if $fg_{i_0}=0$,
we have
$g_{i_0}=0$
since $f$ is a non-zero divisor of $R$.
Thus we have $x_{i_0}g_{i_0}=0$,
this yields a contradiction.

Therefore, we have $fg_{i_0} \neq 0$.
Note that $\deg(fg_{i_0}) < \deg(x_{i_0}g_{i_0})=\indeg_R(\MM_R)$
because $\deg(f) < \deg(x_{i_0})$.
Since $fg_{i_0}$ is a non-zero homogeneous element and $\deg(fg_{i_0}) < \indeg_R(\MM_R)$, we have $fg_{i_0} \in R_0 \setminus \{0\}$.
Thus $fg_{i_0}$ is a unit of $R$ because $R_0$ is a field.
Let $\varphi : I \rightarrow R(-\indeg_R(I)))$ denote the multiplication map by $g_{i_0} \in I^{-1}$.
Then $\varphi$ is a graded $R$-epimorphism of degree $0$ since $\varphi(f) = fg_{i_0}$ is a unit in $R$.  
Let $\psi: R(-\indeg_R(I)) \to I$ be the graded $R$-homomorphism defined by $\psi(1_R) = fg_{i_0}$.
Then $\psi \varphi = \operatorname{id}_I$, and hence $I$ has an $R(-\indeg_R(I))$-free summand.
\end{proof}

The following corollary is a generalization of \cite[Theorem 3.5]{miyashita2024comparing}.

\begin{corollary}\label{Thm:General2}
Let $I \subseteq R$ be a graded ideal of $R$.
Suppose that there exists a non-zero divisor $f \in I_{\indeg_R(I)}$ such that $Rf \neq I$.
If $\tr_R(I)$ contains an $R$-regular sequence $\theta_1,\theta_2 \in R_{\indeg_R(\MM_R)}$,
then $I$ has an $R(-\indeg_R(I))$-free summand.
\begin{proof}
This follows from \autoref{Thm:General} by taking $M=I$.
\end{proof}
\end{corollary}


Recall that the definition of torsion-free element.
Let \( M \) be a graded \( R \)-module.  
An element \( x \in M \) is called \emph{torsion-free} if \( rx = 0 \) implies \( r = 0 \) for all \( r \in R \).  

\begin{theorem}\label{Thm:NonCM1}
Assume that \( R \) is Noetherian
and $\omega_R$ contains a torsion-free homogeneous element of degree \( -a_R \).
Let $I \subseteq R$ be a graded ideal of $R$.
Suppose that
there exists a graded $R$-module $N$ such that
$\omega_R \cong I \oplus N$ as a graded $R$-module (up to degree shift),
$\indeg_R(I)<\indeg_R(N)$
and $[\tr_R(N)]_{\indeg_R(\MM_R)}=0$.
If $\tr_R(\omega_R)$ contains an $R$-regular sequence $\theta_1,\theta_2 \in R_{\indeg_R(\MM_R)}$ and \( \dim_{R_0}([\omega_R]_{-a_R}) = 1 \),
then $R$ is quasi-Gorenstein.
\end{theorem}
\begin{proof}
Since \( \indeg_R(I) < \indeg_R(N) \), we have
$\indeg_R(I \oplus N) = \indeg_R(I)$.
Observe that there exists a non-zero divisor \( f \in I_{\indeg_R(I)} \).

Indeed, by assumption, there exists a homogeneous torsion-free element \( g \in \omega_R \) of degree \( -a_R \),  
and a graded \( R \)-module isomorphism
$\phi: \omega_R \rightarrow I \oplus N$.
Define
$f := \phi(g) \in [I \oplus N]_{\indeg_R(I \oplus N)} = I_{\indeg_R(I)}$.
Then \( f \) is a non-zero divisor, since \( \phi \) is an \( R \)-monomorphism and \( g \) is torsion-free.

We show that \( \tr_R(\omega_R) = R \). Since \( \omega_R \cong I \oplus N \), we have \( \tr_R(\omega_R) = \tr_R(I) + \tr_R(N) \). Thus, it suffices to show that \( \tr_R(I) = R \).

If \( I = Rf \), then \( f \) is a homogeneous non-zero divisor, so \( I \) is isomorphic to \( R \) as a graded \( R \)-module (up to degree shift). In this case, it follows that \( \tr_R(I) = \tr_R(R) = R \).

If \( I \neq Rf \), then \( I \) has an \( R(-\indeg_R(I)) \)-free summand by \autoref{Thm:General}. Thus, we have \( \tr_R(I) = R \).

Therefore, we conclude that \( \tr_R(\omega_R) = R \), and hence \( R \) is quasi-Gorenstein by \cite[Remark 2.9~(4)]{kumashiromiyashita2025}.
\end{proof}

Based on the above preparations, we are now in a position to establish the non-Cohen--Macaulay version of our main result.
Recall that a Noetherian positively graded ring \( R \) is said to be \emph{generically Gorenstein} if \( R_\pp \) is Gorenstein for every \( \pp \in \Ass(R) \); note that we do not assume that \( R \) is Cohen--Macaulay.

\begin{theorem}\label{Thm:NonCM2}
Assume that $R$ is a Noetherian generically Gorenstein graded ring, that $\omega_R$ contains a torsion-free homogeneous element of degree $-a_R$, and that
$\tr_R(\omega_R)$ contains an $R$-regular sequence $\theta_1,\theta_2 \in R_{\indeg_R(\MM_R)}$.
If $\dim_{R_0}([\omega_R]_{-a_R})=1$, then $R$ is quasi-Gorenstein.
\end{theorem}
\begin{proof}
Since \( R \) is generically Gorenstein, it follows from \cite[Remark 2.13~(1)]{kumashiromiyashita2025} that \( \omega_R \) is isomorphic to a graded ideal of \( R \) as a graded \( R \)-module (up to degree shift). Hence, by applying \autoref{Thm:NonCM1}, we conclude that \( R \) is quasi-Gorenstein.
\end{proof}

\begin{remark}
The assumption that \( R \) is a Noetherian generically Gorenstein graded ring such that \( \omega_R \) contains a torsion-free homogeneous element of degree \( -a_R \), as in \autoref{Thm:NonCM2}, is automatically satisfied when \( R \) is a Noetherian graded domain.
\end{remark}

We are now in a position to establish our main result.

\begin{corollary}[{\autoref{XXX}}]\label{mainTHM}
Assume that $R$ is a Cohen--Macaulay generically Gorenstein graded ring, that $\omega_R$ contains a torsion-free homogeneous element of degree $-a_R$, and that
$\tr_R(\omega_R)$ contains an $R$-regular sequence $\theta_1,\theta_2 \in R_{\indeg_R(\MM_R)}$.
If $R$ is pseudo-Gorenstein, then it is Gorenstein.
\end{corollary}
\begin{proof}
This follows from \autoref{Thm:NonCM2}.
\end{proof}

\begin{remark}
\autoref{mainTHM} fails in general if we replace the assumption that $\tr_R(\omega_R)$ contains an $R$-regular sequence $\theta_1, \theta_2 \in R_{\indeg_R(\MM_R)}$ with the weaker condition that there exists a non-zero divisor $\theta \in R_{\indeg_R(\MM_R)}$ such that $\theta \in \tr_R(\omega_R)$.
Indeed, any one-dimensional pseudo-Gorenstein graded domain (e.g., numerical semigroup rings), which is non-Gorenstein and nearly Gorenstein, provides a counterexample.
\end{remark}

\section{Some applications}\label{sectC}
In this section, we present several applications of the results obtained in the previous section.
As a first application of \autoref{Thm:NonCM2}, let us compare several generalizations of the Gorenstein property.
Throughout this section, we retain \autoref{setup1}, and we
further assume that $R$ is Noetherian.



\begin{definition}
Let $\mu(\omega_R)$ denote the number of homogeneous minimal generators of \( \omega_R \) as a graded \( R \)-module.
\end{definition}

\begin{theorem}\label{Thm:WOWOWnonCM}
Assume that
$\omega_R$ contains a torsion-free homogeneous element of degree $-a_R$.
Moreover, we assume that
$\depth(R) \ge 2$ and
$\sqrt{[\tr_R(\omega_R)]_{{\indeg_R(\MM_R) }}R} \supseteq \MM_R.$
Then the following hold:
\begin{itemize}
\item[\rm (1)] If $\dim_{R_0}([\omega_R]_{-a_R})=1$,
then $R$ is quasi-Gorenstein;
\item[\rm (2)] If $\mu(\omega_R)=2$, then $\omega_R$ is generated in a single degree as a graded $R$-module.
\end{itemize}
\end{theorem}
\begin{proof}
First, under the given assumptions, we show that \(R\) is generically Gorenstein.
Since \(\omega_R\) contains a torsion-free homogeneous element, we have \(\Ann_R(\omega_R) = (0)\).  
Therefore, \(\tr_R(\omega_R)\) describes the non-quasi-Gorenstein locus of \(R\) (see \autoref{rem:gradedaoyamagoto}).  
Since \(\sqrt{\tr_R(\omega_R)} \supseteq \sqrt{[\tr_R(\omega_R)]_{\indeg_R(\MM_R)} R} \supseteq \mm_R\),
it follows that $R$ is quasi-Gorenstein on the punctured spectrum.
Therefore, since \(\dim(R) > 0\), we conclude that \(R\) is generically Gorenstein.
Next, we prove statements (1) and (2).

(1):
Suppose that $\dim_{R_0}([\omega_R]_{-a_R})=1$.
Without loss of generality, we may assume that $R_0$ is infinite (see \autoref{HelpMeFromMugentai}).
Note that $\sqrt{[\tr_R(\omega_R)]_{\indeg_R(\MM_R)}R}$
is either \(R\) or \(\MM_R\). In the former case, we have \(\tr_R(\omega_R) = R\), which implies that \(R\) is quasi-Gorenstein~(see \autoref{rem:gradedaoyamagoto}).
In the latter case, since \(\depth\bigl(R) \ge 2\) and \({R_0}\) is an infinite field, we can choose an \(R\)-regular sequence of length two in \([\tr_R(\omega_R)]_{\indeg_R(\MM_R)}\) by the same argument used in the proof of \cite[Proposition 1.5.12]{bruns1998cohen}.
Thus, applying \autoref{Thm:NonCM2} in this situation completes the proof.

(2):
Assume that $\omega_R$ has a minimal generator consisting of two homogeneous elements of distinct degrees as a graded $R$-module.
Then we have $\dim_{R_0}([\omega_R]_{-a_R})=1$, and it follows from \autoref{Thm:WOWOWnonCM}~(1) that \(R\) is quasi-Gorenstein.
This implies that \(\mu(\omega_R)=1<2\), which yields a contradiction.
\end{proof}

\begin{notation}
Suppose that $R$ is Cohen--Macaulay.
Let $r(R)$ denote {\it the Cohen--Macaulay type} of $R$.
\end{notation}

\begin{corollary}\label{Thm:WOWOW}
Assume that $R$ is Cohen--Macaulay,
and that $\omega_R$ contains a torsion-free homogeneous element of degree $-a_R$.
Moreover, we assume that $\dim(R) \ge 2$ and
$\sqrt{[\tr_R(\omega_R)]_{{\indeg_R(\MM_R) }}R} \supseteq \MM_R.$
Then the following hold:
\begin{itemize}
\item[\rm (1)] If $R$ is pseudo-Gorenstein,
then $R$ is Gorenstein;
\item[\rm (2)] If $r(R)=2$, then $R$ is level.
\end{itemize}
\end{corollary}
\begin{proof}
Note that \( r(R)=\mu(\omega_R)\).
Then this follows from \autoref{Thm:WOWOWnonCM}.
\end{proof}

\begin{remark}\label{interestingEXample}
The assertion of \autoref{Thm:WOWOW} does not hold in general if we replace 
$\sqrt{[\tr_R(\omega_R)]_{\indeg_R(\MM_R)} R} \supseteq \MM_R$
with \(\tr_R(\omega_R) \supseteq \mm_R\).
This can be seen, for example, by considering the 2-dimensional affine semigroup ring
$R = \mathbb{Q}[x, xy, x^2y^3, x^3y^5]$,
where the grading on \(R\) is given by
\[
\deg(x) = \deg(xy) = 1, \quad \deg(x^2y^3) = 2, \quad \deg(x^3y^5) = 3.
\]
In this setting, \(R\) is normal (this can be verified using \texttt{Macaulay2}~\cite{M2}).
It is known that a normal affine semigroup ring is Cohen--Macaulay, and that its canonical module is generated by the monomials corresponding to the lattice points in the interior of the rational cone associated to the affine semigroup (see \cite[Theorem 6.3.5]{bruns1998cohen}).
Hence, \(\omega_R \cong I_R = (xy, x^2y^3)R\), and thus \(R\) is a non-level pseudo-Gorenstein ring with $r(R)=2$.
Moreover, since \(\tr_R(\omega_R) = I_R \cdot I_R^{-1} = \mm_R\) (see \cite[Lemma 1.1]{herzog2019trace}), \(R\) is nearly Gorenstein. 
As can be seen from \autoref{Thm:WOWOW},
$\sqrt{[\tr_R(\omega_R)]_{\indeg_R(\MM_R)} R} 
=(x,xy,x^2y^3)R$
does not contain \(\mm_R\).
\end{remark}

On the other hand, the following holds.

\begin{corollary}[{\autoref{YYY}~(1), (2)}]\label{YYYY}
Assume $R$ is Cohen--Macaulay.
Moreover, suppose that \(R\) is nearly Gorenstein with \(\dim(R) \ge 2\), that the ideal generated by \(R_{\indeg_R(\mm_R)}\) is \(\mm_R\)-primary, and that \(\omega_R\) contains a torsion-free homogeneous element of degree \(-a_R\)
(e.g., when \(R\) is a domain).
Then the following hold:
\begin{itemize}
\item[\rm (1)] If $R$ is pseudo-Gorenstein, then it is Gorenstein;
\item[\rm (2)] If $r(R)=2$, then $R$ is level;
\end{itemize}
\begin{proof}
Since \( R \) is nearly Gorenstein and the ideal generated by \( R_{\indeg_R(\mm_R)} \) is \(\mm_R\)-primary, note that
$\sqrt{[\tr_R(\omega_R)]_{\indeg_R(\mm_R)} R} \supseteq \mm_R$.
Then this follows from \autoref{Thm:WOWOW}.
\end{proof}
\end{corollary}

Next, we prove the following statement, which generalizes \cite[Corollary 3.7 and Theorem 6.1]{miyashita2024comparing}.

\begin{corollary}[{\autoref{YYY}~(3)}]\label{Thm:MainC}
Assume that $R$ is Cohen--Macaulay.
Moreover, suppose that $R$ is a semi-standard graded nearly Gorenstein domain with $\dim(R) \ge 2$.
The following hold:
\begin{itemize}
\item[\rm (1)]
If $R$ is pseudo-Gorenstein, then it is Gorenstein;
\item[\rm (2)]
If $r(R)=2$, then $R$ is level;
\item[\rm (3)]
Moreover, assume that \( R \) is almost Gorenstein and standard graded, and \( R_0\) is an algebraically closed field of characteristic zero.
Then \( R \) must be either Gorenstein or have minimal multiplicity.
In particular, $R$ is level in either case.
\end{itemize}
\end{corollary}
\begin{proof}
(1) and (2):
Since \(R\) is semi-standard graded, we have \(\indeg_R(\mm_R) = 1\), and the ideal generated by \(R_1\) is \(\mm_R\)-primary~(see \cite[Chapter I, 5.2]{stanley2007combinatorics}).
Therefore, these follow from \autoref{YYYY}.

(3):
By \autoref{WOW}, \( R \) is either pseudo-Gorenstein or has minimal multiplicity.  
If \( R \) is pseudo-Gorenstein, then it is Gorenstein by \autoref{Thm:MainC}\,(1).  
Therefore, \( R \) is either Gorenstein or has minimal multiplicity.  
It is well known that a Cohen--Macaulay standard graded ring of minimal multiplicity is level.  
For example, this follows from the fact that such a ring has a 2-linear resolution (see \cite[Section 4]{eisenbud1984linear}).
\end{proof}

Finally, we discuss an application to Veronese subalgebras.

\begin{theorem}\label{Thm:MainD}
Assume that \( R \) is a semi-standard graded domain.  
Moreover, suppose that \( a_R \in k\mathbb{Z} \) and \( \depth(R^{(k)}) \geq 2 \) for some \( k \in \mathbb{Z}_{>0} \).  
If \( \operatorname{tr}_R(\omega_R) \) contains an \( \mathfrak{m}_R \)-primary graded ideal $I$ generated by a subset of \( R_1 \), and if \( \dim_{R_0}([\omega_R]_{-a_R}) = 1 \),  
then \( R^{(k)} \) is quasi-Gorenstein.
\begin{proof}
Let $\nn_R$ be the graded ideal generated by $R_1$.
Since \( R \) is semi-standard graded and \( \depth(R) > 0 \), we can take a non-zero divisor \( z \in R_1 \) (see \cite[Lemma 2.18]{miyashita2024linear}).
Then we have
\begin{equation}\label{bbb}
\tr_{R^{(k)}}(\omega_{R ^{(k)}}) \supseteq (\nn_R^{k-1}I)^{(k)}=(\nn_R^{k-1}I) \cap R^{(k)}
\end{equation}
by \autoref{lem:Yeah}.
Since both of $\nn_R$ and $I$ are $\mm_R$-primary ideal~(see \cite[Chapter I, 5.2]{stanley2007combinatorics}), we obtain $\sqrt{\nn_R^{k-1}I}=\mm_R$.
Therefore, we have 
\begin{equation}\label{bbbb}
\sqrt{(\mathfrak{n}_R^{k-1}I) \cap R^{(k)}} = \sqrt{\mathfrak{n}_R^{k-1}I} \cap R^{(k)} = \mathfrak{m}_R \cap R^{(k)} = \mathfrak{m}_{R^{(k)}}.
\end{equation}
Since
\autoref{bbb} and \autoref{bbbb},
and
$(\mathfrak{n}_R^{k-1}I) \cap R^{(k)}$
is generated by a subset of \( [R^{(k)}]_1\), we obtain
$$\sqrt{[\tr_{R^{(k)}}(\omega_{R ^{(k)}})]_1 R^{(k)}} \supseteq 
\sqrt{[(\nn_R^{k-1}I) \cap R^{(k)}]_1 R^{(k)}}=
\sqrt{(\nn_R^{k-1}I) \cap R^{(k)}}=
\mm_{R^{(k)}}.$$
Since \( R \) is a semi-standard graded domain, so is \( R^{(k)} \).

Moreover, from \( \omega_{R^{(k)}}=(\omega_R)^{(k)} \)~(see \cite[Corollary $(3.1.3)$]{goto1978graded}) and \( a_R \in k\mathbb{Z} \), it follows that \( a_{R^{(k)}} = a_R \) and $\dim_{[R^{(k)}]_0}([\omega_{R^{(k)}}]_{-a_{R^{(k)}}}) =
\dim_{R_0}([\omega_R]_{-a_R})=1$.
Therefore, since
$\sqrt{ [\tr_{R^{(k)}}(\omega_{R ^{(k)}})]_1 R^{(k)}} \supseteq \mm_{R^{(k)}}$
and \( \depth(R^{(k)}) \geq 2 \),
we conclude that \( R^{(k)} \) is quasi-Gorenstein by \autoref{Thm:WOWOWnonCM}~(1).
\end{proof}
\end{theorem}

\begin{corollary}[{\autoref{ZZZ}}]\label{OKOKOK}
Let \( R \) be a Noetherian standard graded domain, and let \( k \) be a positive integer.  
Assume that \( a_R \in k\mathbb{Z} \) and that \( \depth(R^{(k)}) \geq 2 \).  
Suppose further that
\[
  \tr_R(\omega_R) \supseteq \MM_R
  \quad \text{and} \quad
  \dim_{R_0}([\omega_R]_{-a_R}) = 1.
\]
Then \( R^{(k)} \) is quasi-Gorenstein.  
In particular, if \( R^{(k)} \) is Cohen--Macaulay, then it is Gorenstein.
\begin{proof}
This follows from \autoref{Thm:MainD}.
\end{proof}
\end{corollary}

\subsection*{Acknowledgement}
I would like to thank Aldo Conca, Kazufumi Eto, Akihiro Higashitani, Shinya Kumashiro, Naoyuki Matsuoka, Takahiro Numata, Kei-ichi Watanabe
and Kohji Yanagawa for their useful comments.
I would also like to express my gratitude to Max K\"olbl for his assistance in proofreading the English manuscript.
This work was supported by JST SPRING, Japan Grant Number JPMJSP2138.



\end{document}